\pdfoutput=1
\RequirePackage{ifpdf}
\ifpdf % We are running pdfTeX in pdf mode
\documentclass[pdftex]{sigma}
\else
\documentclass{sigma}
\fi

\begin{document}

\allowdisplaybreaks

\renewcommand{\PaperNumber}{017}

\FirstPageHeading

\ShortArticleName{Dynamics of an Inverting Tippe Top}

\ArticleName{Dynamics of an Inverting Tippe Top}

\Author{Stefan RAUCH-WOJCIECHOWSKI and Nils RUTSTAM} \AuthorNameForHeading{S.~Rauch-Wojciechowski and
N.~Rutstam} \Address{Department of Mathematics, Link\"{o}ping University, Link\"{o}ping, Sweden}
\Email{\href{mailto:strau@mai.liu.se}{strau@mai.liu.se},
\href{mailto:ergoroff@hotmail.com}{ergoroff@hotmail.com}}

\ArticleDates{Received September 05, 2013, in f\/inal form February 18, 2014; Published online February 27, 2014}

\Abstract{The existing results about inversion of a~tippe top (TT) establish stability of asymptotic
solutions and prove inversion by using the LaSalle theorem.
Dynamical behaviour of inverting solutions has only been explored numerically and with the use of certain
perturbation techniques.
The aim of this paper is to provide analytical arguments \mbox{showing} oscillatory behaviour of TT through the
use of the main equation for the TT.
The main equation describes time evolution of the inclination angle $\theta(t)$ within an ef\/fective
potential $V(\cos\theta,D(t),\lambda)$ that is deforming during the inversion.
We prove here that $V(\cos\theta,D(t),\lambda)$ has only one minimum which (if Jellett's integral is above
a~threshold value
$\lambda>\lambda_{\text{thres}}=\frac{\sqrt{mgR^3I_3\alpha}(1+\alpha)^2}{\sqrt{1+\alpha-\gamma}}$ and
$1-\alpha^2<\gamma=\frac{I_1}{I_3}<1$ holds) moves during the inversion from a~neighbourhood of $\theta=0$
to a~neighbourhood of $\theta=\pi$.
This allows us to conclude that $\theta(t)$ is an oscillatory function.
Estimates for a~maximal value of the oscillation period of $\theta(t)$ are given.}


\Classification{70F40; 74M10; 70E18; 70E40; 37B25}

\section{Introduction}

A tippe top (TT) is constructed as a~truncated axisymmetric sphere with a~small peg as its handle.
The top is spun on a~f\/lat surface with the peg pointing upward.
If the initial rotation is fast enough, the top will start to turn upside down until it ends up spinning on
its peg.
We call this interesting phenomenon an inversion.

It is known that the TT inverts when the physical parameters satisfy the conditions
$1-\alpha<\gamma=\frac{I_1}{I_3}<1+\alpha$ where $0<\alpha<1$ is the eccentricity of the center of mass and
$I_1$, $I_3$ are the main moments of inertia.

The TT and the inversion phenomenon has been studied extensively throughout the years, but the dynamics of
inversion has proven to be a~dif\/f\/icult problem.
This is because even the most simplif\/ied model for the rolling and gliding TT is a~non-integrable
dynamical system with at least 6 degrees of freedom.
The focus in many works has been on the asymptotics of the TT \cite{Mars,Eben,Karap2,Karap5,RSG} or on
numerical simulations for a~TT~\cite{Coh,Or,Ued}.

In this paper we study equations of motion for a~rolling and gliding TT in the case of inverting solutions
and analyse dynamical properties of such solutions through the main equations for the TT~\cite{Rau, Nisse}.

We study the main equation for the TT for a~subset of parameters satisfying $1-\alpha^2<\gamma<1$
and $\frac{1-\gamma}{\gamma+\alpha^2-1}=\frac{mR^2}{I_3}$ when it acquires a~simpler form, which enables
detailed analysis of deformation of the ef\/fective potential $V(\cos\theta,D,\lambda)$ during the
inversion.
We show that, during the inversion, a~minimum of the ef\/fective potential moves from the neighbourhood of
$\theta=0$ to the neighbourhood of $\theta=\pi$ and therefore the inclination angle $\theta(t)$ oscillates
within a~nutational band that moves from the north pole to the south pole of the unit sphere $S^2$.
We give also estimates for the period of nutation of the symmetry axis.

\section{The tippe top model}

We model the TT as an axisymmetric sphere of mass $m$ and radius $R$ which is in instantaneous contact with
the supporting plane at the point $A$.
The center of mass $CM$ is shifted from the geometric center $O$ along its symmetry axis by $\alpha R$,
where $0<\alpha<1$.

We choose a~f\/ixed inertial reference frame $(\widehat{X},\widehat{Y},\widehat{Z})$ with $\widehat{X}$ and
$\widehat{Y}$ parallel to the supporting plane and with vertical $\widehat{Z}$.
We place the origin of this system in the supporting plane.
Let $(\hat{x},\hat{y},\hat{z})$ be a~frame def\/ined through rotation around $\widehat{Z}$ by an angle
$\varphi$, where $\varphi$ is the angle between the plane spanned by $\widehat{X}$ and $\widehat{Z}$ and
the plane spanned by the points $CM$, $O$ and $A$.

\begin{figure}[t]
\centering
\includegraphics[scale=0.90]{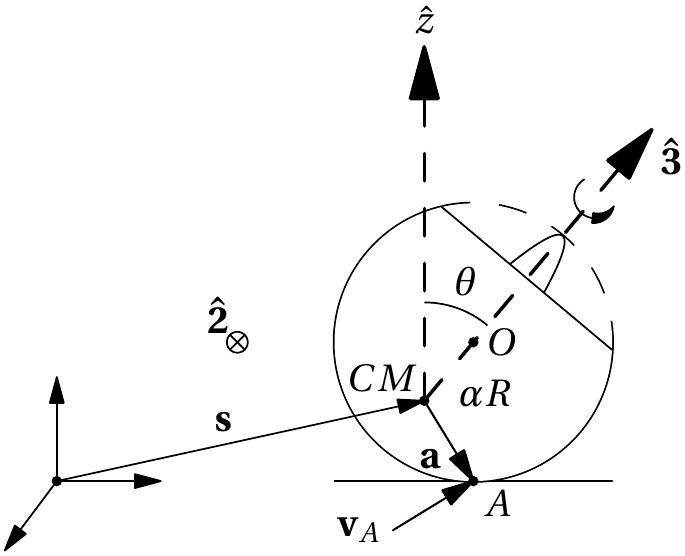}
\caption{Diagram of the TT. Note that $\mathbf{a}=R\alpha\mathbf{\hat{3}}-R\hat{z}$.}\label{TT_diagram_3}
\end{figure}

The third reference frame $(\mathbf{\hat{1}},\mathbf{\hat{2}},\mathbf{\hat{3}})$, with origin at $CM$, is
def\/ined by rotating $(\hat{x},\hat{y},\hat{z})$ by an angle $\theta$ around $\hat{y}$.
Thus $\mathbf{\hat{3}}$ will be parallel to the symmetry axis, and $\theta$ will be the angle between~$\hat{z}$ and~$\mathbf{\hat{3}}$.
This frame is not fully f\/ixed in the body.
The axis $\mathbf{\hat{2}}$ points behind the plane of the picture of Fig.~\ref{TT_diagram_3}.

We let $\mathbf{s}$ denote the position of $CM$ w.r.t.\
the origin of the frame $(\widehat{X},\widehat{Y},\widehat{Z})$ and the vector from $CM$ to $A$ is
$\mathbf{a}=R(\alpha\mathbf{\hat{3}}-\hat{z})$.
The orientation of the body w.r.t.\
the inertial reference frame $(\hat{X},\hat{Y},\hat{Z})$ is described by the Euler angles
$(\theta,\varphi,\psi)$, where $\psi$ is the rotation angle of the sphere about the symmetry axis.
With this notation, the angular velocity of the TT is
$\boldsymbol{\omega}=-\dot{\varphi}\sin\theta\mathbf{\hat{1}}+\dot{\theta}\mathbf{\hat{2}}+(\dot{\psi}+\dot{\varphi}\cos\theta)\mathbf{\hat{3}}$,
and we denote $\omega_3:=\dot{\psi}+\dot{\varphi}\cos\theta$.

The principal moments of inertia along the axes $(\mathbf{\hat{1}},\mathbf{\hat{2}},\mathbf{\hat{3}})$ are
denoted by $I_1=I_2$ and $I_3$, so the inertia tensor $\mathbb{I}$ will have components $(I_1,I_1,I_3)$
with respect to the $(\mathbf{\hat{1}},\mathbf{\hat{2}},\mathbf{\hat{3}})$-frame.
The axes $\mathbf{\hat{1}}$ and $\mathbf{\hat{2}}$ are principal axes due to the axisymmetry of TT.
The equations of motion for TT are the Newton equations for the rolling and gliding rigid body
\begin{gather}
\label{TT_equ}
m\mathbf{\ddot{s}}=\mathbf{F}-mg\hat{z},
\qquad
\mathbf{\dot{L}}=\mathbf{a\times F},
\qquad
\mathbf{\dot{\mathbf{\hat{3}}}}=\boldsymbol{\omega}\times\mathbf{\hat{3}},
\end{gather}
where $\mathbf{F}$ is the external force acting on the TT at the supporting point $A$ and
$\mathbf{L}=\mathbb{I}\boldsymbol{\omega}$ is the angular momentum w.r.t.\
 $CM$.
We assume that the TT is always in contact with the plane at $A$, so $\hat{z}\cdot(\mathbf{a+s})=0$ holds
at all times.
This system is known to admit Jellett's integral of motion ${\lambda=-\mathbf{L\cdot
a}=RI_1\dot{\varphi}\sin^2\theta-RI_3\omega_3(\alpha-\cos\theta)}$ (without loss of generality, we will
assume in this paper that $\lambda$ is positive).

The contact condition determines the vertical part of the external force, but the planar parts must be
specif\/ied to make system~\eqref{TT_equ} complete.
We assume that the contact force has the form $\mathbf{F}=g_n\hat{z}-\mu g_n\mathbf{v}_A$, where $g_n\geq
0$ is the normal force and $-\mu g_n\mathbf{v}_A$ is a~viscous-type friction force, acting against the
gliding velocity $\mathbf{v}_A$.
The quantity $\mu(\mathbf{L},\mathbf{\hat{3}},\mathbf{\dot{s}},\mathbf{s},t)\geq 0$ is a~friction
coef\/f\/icient.

For this model of the rolling and gliding TT, it is easy to see~\cite{DelCampo,Nisse2} that the energy is
decreasing $\dot{E}=\mathbf{F}\cdot\mathbf{v}_A<0$ and that the $\hat{y}$ component of the friction force
is the only force creating the torque necessary for transferring the rotational energy into the potential
energy, thus lifting the $CM$ of the TT.
This mechanism shows that the inversion phenomenon is created by the gliding friction.

The asymptotic properties of this model have been analysed in previous
works~\cite{Mars,Eben,Karap2,Karap5,RSG}.
In the nongliding case, $\mathbf{v}_A=0$, the possible motions for the TT are either spinning in the
upright ($\theta=0$) or in the inverted ($\theta=\pi$) position, or rolling around with f\/ixed $CM$ with
an inclination angle $\theta\in(0,\pi)$.

The inclined rolling solutions are called tumbling solutions.
If ${1-\alpha<\gamma<1+\alpha}$, where $\gamma=I_1/I_3$, every angle in the interval $(0,\pi)$ determines
an admissible tumbling solution.
Further, by a~LaSalle-type theorem~\cite{RSG}, it is known that for initial conditions such that the
absolute value of the Jellett integral $|\lambda|$ is above the threshold value $\lambda_{\text{thres}}$,
only the inverted spinning position is a~stable asymptotic solution.
For a~TT built such that it satisf\/ies the parameter condition $1-\alpha<\gamma<1+\alpha$ and for initial
conditions with $\mathbf{L}\cdot\hat{z}$ such that ${\lambda>\lambda_{\text{thres}}}$, the inversion can
take place.

Since we are primarily interested in the dynamics of inversion and we want to consider solutions describing
an inverting TT, the basic assumptions are that the TT in question satisf\/ies the parameter constraint and
that we have initial conditions such that $\lambda$ is above the \mbox{threshold}.
Then we have a~situation where an inverting solution becomes the only stable asymptotic solution, so the TT
has to invert.
Our aim is to describe the dynamics of inverting solutions.

In our particular model, the assumptions about the reaction force $\mathbf{F}$ and the contact constraint
yield the reduced equations of motion for the rolling and gliding TT:
\begin{gather}
\label{TT_equ_red}
\frac{d}{dt}\left(\mathbb{I}\boldsymbol{\omega}\right)=\mathbf{a}\times\left(g_n\hat{z}-\mu g_n\mathbf{v}
_A\right),
\qquad
m\mathbf{\ddot{r}}=-\mu g_n\mathbf{v}_A,
\qquad
\mathbf{\dot{\mathbf{\hat{3}}}}=\boldsymbol{\omega}\times\mathbf{\hat{3}},
\end{gather}
where $\mathbf{r}=\mathbf{s}-s_{\hat{z}}\hat{z}$.
We write the gliding velocity as
$\mathbf{v}_A=\nu_{x}\cos\theta\mathbf{\hat{1}}+\nu_y\mathbf{\hat{2}}+\nu_x\sin\theta\mathbf{\hat{3}}$,
where~$\nu_x$,~$\nu_y$ are the velocities in the $\mathbf{\hat{2}}\times\hat{z}$ and $\mathbf{\hat{2}}$
direction.
Equations~\eqref{TT_equ_red} can be written in the Euler form and then solved for the highest derivative of
each of the variables $(\theta,\varphi,\omega_3,\nu_x,\nu_y)$.
We then get the system
\begin{gather}
\label{ddth3}
\ddot{\theta}=\frac{\sin\theta}{I_1}\left(I_1\dot{\varphi}^2\cos\theta-I_3\omega_3\dot{\varphi}
-R\alpha g_n\right)+\frac{R\mu g_n\nu_x}{I_1}(1-\alpha\cos\theta),
\vspace{-3pt}\\
\label{ddph3}
\ddot{\varphi}=\frac{I_3\dot{\theta}\omega_3-2I_1\dot{\theta}\dot{\varphi}
\cos\theta-\mu g_n\nu_y R(\alpha-\cos\theta)}{I_1\sin\theta},
\\
\label{dot_om3}
\dot{\omega}_3=-\frac{\mu g_n\nu_y R\sin\theta}{I_3},
\\
\dot{\nu}_x=\frac{R\sin\theta}{I_1}\left(\dot{\varphi}
\omega_3\left(I_3(1-\alpha\cos\theta)-I_1\right)+g_nR\alpha(1-\alpha\cos\theta)-I_1\alpha\big(\dot{\theta}
^2+\dot{\varphi}^2\sin^2\theta\big)\right)
\nonumber
\\
\phantom{\dot{\nu}_x=}
{}-\frac{\mu g_n\nu_x}{mI_1}\left(I_1+mR^2(1-\alpha\cos\theta)^2\right)+\dot{\varphi}\nu_y,
\label{nu_x3}
\\
\dot{\nu}_y=-\frac{\mu g_n\nu_y}{mI_1I_3}
\left(I_1I_3+mR^2I_3(\alpha-\cos\theta)^2+mR^2I_1\sin^2\theta\right)
\nonumber
\\
\phantom{\dot{\nu}_y=}
{}+\frac{\omega_3\dot{\theta}R}{I_1}\left(I_3(\alpha-\cos\theta)+I_1\cos\theta\right)-\dot{\varphi}\nu_x,
\label{nu_y3}
\end{gather}
which, if we add the equation $\frac{d}{dt}(\theta)=\dot{\theta}$, becomes a~dynamical system of the form
%$
\begin{gather*}
(\dot{\theta},\ddot{\theta},\ddot{\varphi},\dot{\omega}_3,\dot{\nu}_x,\dot{\nu}_y)=(h_1(\theta,\ldots,\nu_y),\ldots,h_6(\theta,\ldots,\nu_y)).
\end{gather*}
The value of the normal force $g_n$ can be determined from the contact constraint
$(\mathbf{a+s})\cdot\hat{z}=0$ and it is{\samepage
\begin{gather}
\label{equ_for_g_n}
g_n=\frac{mgI_1+mR\alpha(\cos\theta(I_1\dot{\varphi}^2\sin^2\theta+I_1\dot{\theta}^2)-I_3\dot{\varphi}
\omega_3\sin^2\theta)}{I_1+mR^2\alpha^2\sin^2\theta-mR^2\alpha\sin\theta(1-\alpha\cos\theta)\mu\nu_x}.
\end{gather}
We see that we get a~complicated, nonlinear system for 6 unknowns.}

\section{The main equation for the tippe top}

For further study of inverting solutions we need to clarify the logic of applying the main equation for the
tippe top (METT) to analysing motion of TT.
We need also to recall properties of TT equations when the TT is only rolling on the supporting surface and
the gliding velocity vanishes $\mathbf{v}_A=\mathbf{\dot{s}}+\boldsymbol{\omega}\times\mathbf{a}=0$.
It is the well known~\cite{Chap2,Glad} integrable case of the rolling axisymmetric sphere that was f\/irst
separated by Chaplygin.
We need to explain how the structure of separation equations motivates the introduction of the METT and how
this equation dif\/fers from the classical separation equation.

For the purely rolling axisymmetric sphere, the constraint
$\mathbf{v}_A=\mathbf{\dot{s}}+\boldsymbol{\omega}\times\mathbf{a}=0$ implies that the equations of
motion~\eqref{TT_equ} reduce to a~closed system for the vectors $\mathbf{\hat{3}}$ and
$\boldsymbol{\omega}$:
\begin{gather}
\label{rTT_equ}
\frac{d}{dt}\left(\mathbb{I}\boldsymbol{\omega}\right)=m\mathbf{a}\times\left(g\hat{z}-\frac{d}{dt}
(\boldsymbol{\omega}\times\mathbf{a})\right),
\qquad
\mathbf{\dot{\mathbf{\hat{3}}}}=\boldsymbol{\omega}\times\mathbf{\hat{3}}.
\end{gather}
For this system the external force is dynamically determined:
${\mathbf{F}=mg\hat{z}-m\frac{d}{dt}(\boldsymbol{\omega}\times\mathbf{a})}$.
In the Euler angle form the equations give a~fourth order dynamical system for
$(\theta,\dot{\theta},\dot{\varphi},\omega_3)$.

The system~\eqref{rTT_equ} admits three integrals of motion.
Since the system is conservative, the energy
\begin{gather}
E=\frac{1}{2}m\mathbf{\dot{s}}^2+\frac{1}{2}\boldsymbol{\omega}\cdot\mathbf{L}+mg\mathbf{s}\cdot\hat{z}
=\frac{1}{2}\left(I_1\dot{\varphi}^2\sin^2\theta+I_1\dot{\theta}^2+I_3\omega_{3}
^2\right)+mgR(1-\alpha\cos\theta)
\nonumber
\\
\phantom{E=}
{}+\frac{1}{2}mR^2\bigg[(\alpha-\cos\theta)^2(\dot{\theta}^2+\dot{\varphi}
^2\sin^2\theta)+\sin^2\theta(\dot{\theta}^2+\omega_{3}^2+2\omega_3\dot{\varphi}(\alpha-\cos\theta))\bigg]
\label{Energy3}
\end{gather}
is an integral of motion.
We also have Jellett's integral ${\lambda=RI_1\dot{\varphi}\sin^2\theta-RI_3\omega_3(\alpha-\cos\theta)}$
as well as the Routh integral
\begin{gather*}
D:=\omega_3\sqrt{I_3I_1+mR^2(\alpha-\cos\theta)^2+mR^2I_1\sin^2\theta}=I_3\omega_3\sqrt{d(\cos\theta)},
\end{gather*}
where
\begin{gather*}
d(z)=\gamma+\sigma(\alpha-z)^2+\sigma\gamma(1-z^2),
\qquad
\sigma=\frac{mR^2}{I_3}
\qquad
\text{and}
\qquad
\gamma=\frac{I_1}{I_3}.
\end{gather*}
They allow to eliminate $\omega_3=\frac{D}{I_3\sqrt{d(\cos\theta)}}$ and
$\dot{\varphi}=\frac{\lambda\sqrt{d(\cos\theta)}+RD(\alpha-\cos\theta)}{RI_1\sin^2\theta\sqrt{d(\cos\theta)}}$
from the expression of the energy~\eqref{Energy3} to get the separation equation
\begin{gather}
\label{MErTT3}
E=g(\cos\theta)\dot{\theta}^2+V(\cos\theta,D,\lambda),
\end{gather}
where $g(\cos\theta)=\frac{1}{2}I_3\left(\sigma((\alpha-\cos\theta)^2+1-\cos^2\theta)+\gamma\right)$ and
\begin{gather*}
V(z=\cos\theta,D,\lambda)=mgR(1-\alpha z)+\frac{(\lambda\sqrt{d(z)}+RD(\alpha-z))^2}{2I_3R^2\gamma^2(1-z^2)}
+\frac{(R^2D^2-\sigma\lambda^2)}{2R^2I_1}.
\end{gather*}

The separable f\/irst order dif\/ferential equation~\eqref{MErTT3} for $\theta$ determines the motion of
the rolling TT.
It is the Chaplygin separation equation for an axisymmetric sphere~\cite{Chap2}.
We shall show that for certain choice of parameters the ef\/fective potential $V(z,D,\lambda)$ is convex in
$z\in[-1,1]$ so, since $V(z,D,\lambda)\to\infty$ as $z\to\pm1$, it has one minimum in the interval $[-1,1]$.
This means that for f\/ixed $E$ the solutions $\theta(t)$ describe nutational motion of the rolling TT
between two bounding angles $\theta_1$, $\theta_2$ determined by the equation $E=V(\cos\theta,D,\lambda)$.

The rolling and gliding TT only has $\lambda$ as an integral of motion.
It is useful however to consider $D(\theta(t),\omega_3(t))=I_3\omega_3(t)\sqrt{d(\cos\theta(t))}$ being now
a~time dependent function.
Its derivative we calculate using the equations of motion~\eqref{ddth3}--\eqref{nu_y3} for the rolling and
gliding TT:
\begin{gather*}
\frac{d}{dt}D(\theta,\omega_3)=\frac{\gamma mR\sin\theta}{\sqrt{d(\cos\theta)}}(\dot{\varphi}\nu_x+\dot{\nu}_y)
=\frac{\gamma m}{\alpha\sqrt{d(\hat{z}\cdot\mathbf{\hat{3}})}}(\hat{z}\times\mathbf{a})\cdot\mathbf{\dot{v}}_A.
\end{gather*}
For the total energy of TT,
$E=\frac{1}{2}m(\mathbf{v}_A-\boldsymbol{\omega}\times\mathbf{a})^2+\frac{1}{2}\boldsymbol{\omega}\cdot\mathbf{L}+mg\mathbf{s}\cdot\hat{z}$,
we know that $\dot{E}=\mathbf{F}\cdot\mathbf{v}_A<0$.
The part of $E$ that does not depend on $\mathbf{v}_A$,
\begin{gather}
\tilde{E}(\theta,\dot{\theta},\dot{\varphi},\omega_3)=\frac{1}{2}m(\boldsymbol{\omega}\times\mathbf{a})^2
+\frac{1}{2}\boldsymbol{\omega}\cdot\mathbf{L}+mg\mathbf{s}\cdot\hat{z}
\nonumber
\\
\phantom{\tilde{E}(\theta,\dot{\theta},\dot{\varphi},\omega_3)}
=\frac{1}{2}\left(I_1\dot{\varphi}^2\sin^2\theta+I_1\dot{\theta}^2+I_3\omega_{3}
^2\right)+mgR(1-\alpha\cos\theta)
\label{modified_energy3}
\\
\phantom{\tilde{E}(\theta,\dot{\theta},\dot{\varphi},\omega_3)=}
{}+\frac{1}{2}mR^2\big[(\alpha-\cos\theta)^2(\dot{\theta}^2+\dot{\varphi}
^2\sin^2\theta)+\sin^2\theta(\dot{\theta}^2+\omega_{3}^2+2\omega_3\dot{\varphi}(\alpha-\cos\theta))\big],
\nonumber
\end{gather}
we will call the modif\/ied energy function.
The derivative of this function is
$\frac{d}{dt}\tilde{E}(\theta,\dot{\theta},\dot{\varphi},\omega_3)=m\mathbf{\dot{v}}_A\cdot(\boldsymbol{\omega}\times\mathbf{a})$.

With the use of the functions $D(\theta,\omega_3)$, $\tilde{E}(\theta,\dot{\theta},\dot{\varphi},\omega_3)$
we can write the TT equations of motion~\eqref{ddth3}--\eqref{nu_y3} in an equivalent integrated
form~\cite{Rau,Nisse2} as
\begin{gather*}
\frac{d}{dt}\lambda(\theta,\dot{\theta},\dot{\varphi},\omega_3)=0,
\\
\frac{d}{dt}D(\theta,\omega_3)=\frac{\gamma m}{\alpha\sqrt{d(\hat{z}\cdot\mathbf{\hat{3}})}}(\hat{z}\times\mathbf{a})\cdot\mathbf{\dot{v}}_A,
\\
\frac{d}{dt}\tilde{E}(\theta,\dot{\theta},\dot{\varphi},\omega_3)=m(\boldsymbol{\omega}\times\mathbf{a})\cdot\mathbf{\dot{v}}_A,
\\
\frac{d}{dt}m\mathbf{\dot{r}}=-\mu g_n\mathbf{v}_A.
\end{gather*}
These equations are as dif\/f\/icult as the equations~\eqref{ddth3}--\eqref{nu_y3}.
However, if we treat $D(\theta(t),\omega_3(t))=:D(t)$,
$\tilde{E}(\theta(t),\dot{\theta}(t),\dot{\varphi}(t),\omega_3(t))=:\tilde{E}(t)$ as given known functions,
then from $D(t)=I_3\omega_3\sqrt{d(\cos\theta)}$ and
$\lambda=RI_1\dot{\varphi}\sin^2\theta-RI_3\omega_3(\alpha-\cos\theta)$ we can calculate $\dot{\varphi}$,
$\omega_3$ and substitute into expression~\eqref{modified_energy3} for the modif\/ied energy to obtain the
METT~\cite{Rau,Nisse2} that involves only the function $\theta(t)$:
\begin{gather*}%\label{METT3}
\tilde{E}(t)=g(\cos\theta)\dot{\theta}^2+V(\cos\theta,D(t),\lambda).
\end{gather*}
This equation has the same form as equation~\eqref{MErTT3}, but now it depends explicitly on time through
the functions $D(t)$ and $\tilde{E}(t)$.
Solving this equation is therefore not longer possible.
It is a~f\/irst order time dependent ODE which we can study provided that we have some quantitative
information about the functions $D(t)$ and $\tilde{E}(t)$.
The functions $D(t)$, $\tilde{E}(t)$ are usually unknown but for inverting solutions we have qualitative
information about their behaviour due to conservation of the Jellett function $\lambda$.

Thus we consider the motion of the TT as being determined by the three functions
$(\lambda,D(t)$, $\tilde{E}(t))$ and governed by the METT.

Of particular interest regarding the inversion movement is the initial and f\/inal position of the TT.
The TT goes (asymptotically) from an initial angle close to $\theta=0$ to the f\/inal angle close to
$\theta=\pi$ which means, since $\lambda=-\mathbf{L}\cdot\mathbf{a}$ is constant, that
$\lambda=L_0R(1-\alpha)=L_1R(1+\alpha)$ (where $L_0$ and $L_1$ are the values of $|\mathbf{L}|$ at
$\theta=0$ and $\theta=\pi$, respectively).
This implies that $D_0=L_0\sqrt{d(1)}=\frac{\lambda}{R(1-\alpha)}\sqrt{\gamma+\sigma(1-\alpha)^2}$ and
$D_1=-L_1\sqrt{d(-1)}=-\frac{\lambda}{R(1+\alpha)}\sqrt{\gamma+\sigma(1+\alpha)^2}$, and also that
$\tilde{E}_0=\frac{\lambda^2}{2R^2I_3(1-\alpha)^2}+mgR(1-\alpha)$ and
$\tilde{E}_1=\frac{\lambda^2}{2R^2I_3(1+\alpha)^2}+mgR(1+\alpha)$~\cite{Nisse}.

The values $(D_0,\tilde{E}_0)$ and $(D_1,\tilde{E}_1)$ can be interpreted as the boundary values for the
unknown functions $(D(t),\tilde{E}(t))$.
So we assume that for inverting solutions
$(D(t),\tilde{E}(t))\stackrel{t\to\infty}{\longrightarrow}(D_1,\tilde{E}_1)$ and
$(D(t),\tilde{E}(t))\stackrel{t\to-\infty}{\longrightarrow}(D_0,\tilde{E}_0)$.

The aim of the following sections is to analyse dynamical properties of the inverting solution as the
symmetry axis of TT moves from a~neighborhood of $\theta=0$ to a~neighborhood of $\theta=\pi$.

In order to simplify the technical side of analysis we choose special values of the parameters in METT so
that the ef\/fective potential $V(\cos\theta,D,\lambda)$ becomes rational, but we expect that the whole
line of reasoning can be repeated in the general case when the potential depends algebraically on $z$
through $\sqrt{d(z)}$.

We show that $V(z,D,\lambda)$ is strictly convex and therefore has one minimum $z_{\min}$.
We show also that, for inverting solutions when $D(t)$ moves from $D_0$ to $D_1$, the potential deforms so
that $z_{\min}=z_{\min}(D,\lambda)$ moves from a~neighborhood of $z=1$ to a~neighborhood of $z=-1$.

Thus as the potential $V(z,D(t),\lambda)$ deforms and the modif\/ied energy $\tilde{E}(t)$ changes from
$\tilde{E}_0$ to $\tilde{E}_1$ the angle $\theta(t)$ performs oscillatory motions between two turning
angles $\theta_{\pm}(t)$ satisfying the equation $V(\cos\theta_{\pm},D(t),\lambda)=\tilde{E}(t)$.

On the unit sphere $S^2$ the angle $\theta(t)$ performs nutational motion within the nutational band
$[\theta_{-}(t),\theta_{+}(t)]$ that moves from the neighborhood of the north pole to the neighborhood of
the south pole.

We shall give an estimate for the relation between the inversion time $T_{\text{inv}}$ and the maximal
period of nutation $T_V(\tilde{E}(t),D(t))$, so that if $T_{\text{inv}}$ is an order of magnitude larger
than $T_{V}$, say $T_{\text{inv}}>10T_V$, the angle $\theta(t)$ performs oscillatory motion within the
moving nutational band. %\text

\subsection{The rational form of the METT}

The ef\/fective potential in the separation equation~\eqref{MErTT3} is an algebraic function in $z$, which
complicates the analysis.
We can however make a~restriction on the parameters so that the second degree polynomial $d(z)$ can be
written as a~perfect square.
This makes the term $\sqrt{d(z)}$ a~linear function of $z$ and the potential becomes a~rational
function~\cite{Chap2,Nisse3}.

We see that if $1-\alpha^2<\gamma<1$, and if we let the parameter
$\sigma=\frac{1-\gamma}{\gamma+\alpha^2-1}>0$, then
$d(z)=\gamma+\sigma(\alpha-z)^2+\sigma\gamma(1-z^2)=\frac{(\alpha-(1-\gamma)z)^2}{\gamma+\alpha^2-1}$.
This is a~perfect square, so for $\gamma$ in this range we can f\/ind physical values for $\sigma$ such
that $\sqrt{d(z)}$ is a~real polynomial in $z$: $\frac{\alpha-(1-\gamma)z}{\sqrt{\gamma+\alpha^2-1}}$.
Note that $(1-\alpha^2,1)$ is a~subinterval of $(1-\alpha,1+\alpha)$, the parameter range for $\gamma$
where complete inversion of TT is possible.

When $\sigma=\frac{1-\gamma}{\gamma+\alpha^2-1}$, we can rewrite the functions in the separation equation
$E=g(\cos\theta)\dot{\theta}^2+V(\cos\theta,D,\lambda)$ as
\begin{gather*}
g(z)=\frac{I_3}{2}\big(\sigma\big(1+\alpha^2-2\alpha z\big)+\gamma\big)=\frac{I_3}{2}\frac{1}{\gamma+\alpha^2-1}
\left(\alpha^2+(1-\gamma)^2-2\alpha(1-\gamma)z)\right),
\end{gather*}
and
\begin{gather}
V(z,D,\lambda)=mgR(1-\alpha z)+\frac{\big(\lambda(\alpha-(1-\gamma)z)+RD\sqrt{\gamma+\alpha^2-1}(\alpha-z)\big)^2}
{2I_3R^2\gamma^2(\gamma+\alpha^2-1)(1-z^2)}
\nonumber
\\
\phantom{V(z,D,\lambda)=}
{}+\frac{R^2D^2(\gamma+\alpha^2-1)-(1-\gamma)\lambda^2}{2R^2I_1(\gamma+\alpha^2-1)}.
\label{rational_pot}
\end{gather}
This rational form of the ef\/fective potential is simpler to work with.

We should note that the restriction on the parameter $\sigma$ implies that the moments of inertia~$I_1$ and
$I_3$ are dependent on each other; if $\sigma=\frac{1-\gamma}{\gamma+\alpha^2-1}$ then
$I_1=\frac{I_{3}^2+mR^2I_3(1-\alpha^2)}{I_3+mR^2}<I_3$.

\section[Convexity of the rational potential $V(z,D,\lambda)$]{Convexity of the rational potential $\boldsymbol{V(z,D,\lambda)}$}

The range of parameters making the ef\/fective potential a~rational function provides a~simplest
non-trivial situation in which we can study properties of the potential in greater detail.

\begin{figure}[t]
\centering
\includegraphics[scale=1]{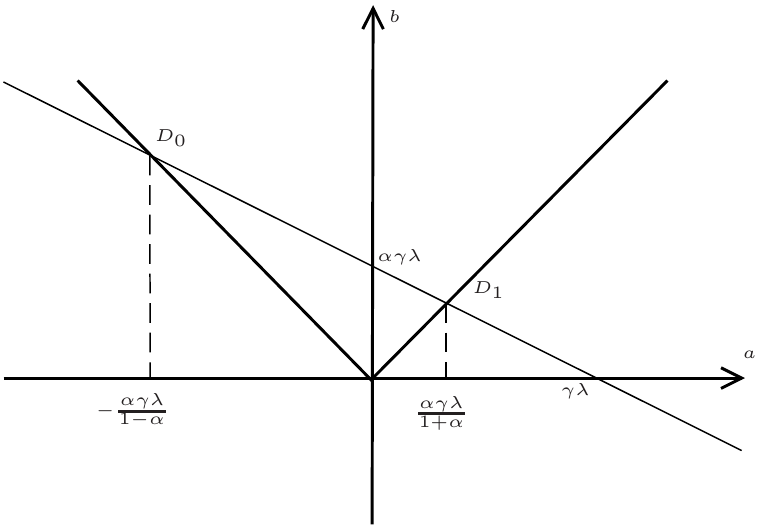}
\caption{Linear relationship between the parameters $a$ and $b$.
This line intersects the lines $b=a$ and $b=-a$ when $D=D_1$ and $D=D_0$, respectively.}\label{ab_diagram}
\end{figure}

We consider the potential function
\begin{gather}
\tilde{V}(z,D,\lambda)=-mgR\alpha z
+\frac{\big(\lambda(\alpha-(1-\gamma)z)+RD\sqrt{\gamma+\alpha^2-1}(\alpha-z)\big)^2}{2I_3R^2\gamma^2(\gamma+\alpha^2-1)(1-z^2)}
\nonumber
\\
\phantom{\tilde{V}(z,D,\lambda)}
=:\frac{1}{2I_3R^2\gamma^2(\gamma+\alpha^2-1)}\left(-\beta z+\frac{(az+b)^2}{1-z^2}\right),
\label{defn_of_fz}
\end{gather}
where
$\tilde{V}(z,D,\lambda)=V(z,D,\lambda)-\big(mgR+\frac{R^2D^2(\gamma+\alpha^2-1)-(1-\gamma)\lambda^2}{2R^2I_1(\gamma+\alpha^2-1)}\big)$,
since the constant does not af\/fect the shape of $V(z,D,\lambda)$ and the position of minimum $z_{\min}$.
The parameters in the function $f(z)=-\beta z+\frac{(az+b)^2}{1-z^2}$ (the expression inside the
parentheses on the r.h.s.\
of~\eqref{defn_of_fz}) are therefore def\/ined as
\begin{gather}
\label{a_def}
a=-(1-\gamma)\lambda-RD\sqrt{\gamma+\alpha^2-1},
\\
\label{b_def}
b=\alpha\lambda+\alpha RD\sqrt{\gamma+\alpha^2-1},
\\
\label{beta_def}
\beta=2mgR^3\alpha I_3\gamma^2(\gamma+\alpha^2-1).
\end{gather}
Remember that $1-\alpha^2<\gamma<1$ and the range of parameters $a$, $b$ is determined by the range of $D$.
We observe that $a+b=0$ if
$D=D_0=\frac{\lambda}{R(1-\alpha)}\frac{(\alpha-1+\gamma)}{\sqrt{\gamma+\alpha^2-1}}$ and $a-b=0$ if
$D=D_1=-\frac{\lambda}{R(1+\alpha)}\frac{(\alpha+1-\gamma)}{\sqrt{\gamma+\alpha^2-1}}$.
The parameters $a$ and $b$ satisfy the relation $b+\alpha a=\lambda\gamma\alpha$.
It is illustrated in Fig.~\ref{ab_diagram},
where the lines $b=a$ and $b=-a$ correspond to $D=D_1$ and $D=D_0$, respectively.

\begin{proposition}\label{rational_convex}
The effective potential $V(z,D,\lambda)$ of~\eqref{rational_pot} is convex for $z\in(-1,1)$ and for all
real values of $D$ and $\lambda$.
\end{proposition}
\begin{proof}
We must show that $\frac{d^2}{dz^2}V(z,D,\lambda)\geq 0$ for $z\in(-1,1)$.
Due to the form of $V$, it is enough to show that the rational function $\frac{(az+b)^2}{1-z^2}$ is convex
for all $a$, $b$.
Suppose f\/irst $a\neq\pm b$ and $ab\neq0$.
We look at the second derivative of this function:
\begin{gather*}
\frac{d^2}{dz^2}\frac{(az+b)^2}{1-z^2}=\frac{2}{(1-z^2)^3}(2abz^3+3(a^2+b^2)z^2+6abz+a^2+b^2)=\frac{2q(z)}{(1-z^2)^3},
\end{gather*}
and have to show that the third degree polynomial $q(z)=2abz^3+3(a^2+b^2)z^2+6abz+a^2+b^2$, has no roots in
the interval $[-1,1]$.

To do this we apply the Sturm theorem~\cite{sturmund}.
We generate a~sequence of polynomials $(q_0(z),q_1(z),q_2(z),\ldots,q_m(z))$ recursively by starting from
a~square-free polynomial $q(z)$: $q_0(z)=q(z)$, $q_1(z)=q'(z)$ and $q_{i}=-\text{rem}(q_{i-1},q_{i-2})$ for $i\geq 2$. %\text
Here $\text{rem}(q_{i-1},q_{i-2})$ denotes the remainder after polynomial division of $q_{i-1}$ by $q_{i-2}$.
By Euclid's algorithm, this will terminate with the constant polynomial $q_m$.
Let $S(\xi)$ be the number of sign-changes in the sequence $(q_0(\xi)$, $q_1(\xi), q_2(\xi),\ldots,q_m(\xi))$
at the point $\xi$.
The Sturm theorem states that for real numbers $c<d$ the number of distinct roots in $(c,d]$ is $S(c)-S(d)$.

For our polynomial $q(z)$ the algorithm described above yields four polynomials $(q_0,q_1,q_2,q_3)$.
When we look at this sequence of polynomials at the points $z=-1$ and $z=1$ we have
\begin{alignat*}{3}
& q_0(-1)=4(a-b)^2,
\qquad &&
q_0(1)=4(a+b)^2, &
\\
& q_1(-1)=-6(a-b)^2,
\qquad &&
q_1(1)=6(a+b)^2,&
\\
& q_2(-1)=-\frac{(a^2-b^2)^2}{ab},
\qquad &&
q_2(1)=\frac{(a^2-b^2)^2}{ab}, &
\\
& q_3(-1)=-6ab,
\qquad &&
q_{3}(1)=-6ab. &
\end{alignat*}
We thus see that the number of sign changes for the Sturm sequence at both points $-1$ and $1$ is the same,
either 1 or 2, depending on whether $ab$ is positive or negative.
Thus according to Sturm's theorem, $q(z)$ has no roots in $(-1,1]$.
Since $q(0)=a^2+b^2>0$ (and $q(-1)>0$), $q(z)$ is positive in $[-1,1]$ and we can conclude that
$\frac{(az+b)^2}{1-z^2}$ is convex when $a\neq\pm b$ and $ab\neq0$.
Routine checking conf\/irms that this is also convex if $a=\pm b$ and as well as if $ab=0$.
Thus the function is convex for all arbitrary values of the parameters $a,\;b\in\mathbb{R}$.
\end{proof}

\subsection[Estimates for position of minimum of $V(z,D,\lambda)$]{Estimates for position of minimum of
$\boldsymbol{V(z,D,\lambda)}$}

The result that $V(z,D,\lambda)$ is convex for $z\in(-1,1)$, combined with ${V(z,D,\lambda)\to\infty}$ as
$z\to\pm 1$, gives that the potential will have exactly one minimum $z_{\min}$ in $[-1,1]$.
We use this to analyse how this minimum shifts as $D$ goes from $D_0$ to $D_1$, for $\lambda>\lambda_{\text{thres}}$.
We show below that for every (small) $\epsilon>0$ we can f\/ind $\delta_{+}(\epsilon,\lambda)>0$ and
$\delta_{-}(\epsilon,\lambda)>0$ such that for $0<\delta<\delta_{+}(\epsilon,\lambda)$ and
$0<\delta<\delta_{-}(\epsilon,\lambda)$ the potential
$V(z,D_0\mp\delta/(R(1-\alpha)\sqrt{\gamma+\alpha^2-1}),\lambda)$ has a~minimum in $[1-\epsilon,1]$ and
$V(z,D_1\pm\delta/(R(1+\alpha)\sqrt{\gamma+\alpha^2-1}),\lambda)$ has a~minimum in $[-1,-1+\epsilon]$.

This means that as the inclination angle $\theta$ goes from a~value close to $0$ to a~value close to $\pi$
and $D(t)$ goes from a~neighborhood
$\left(D_0-\frac{\delta}{R(1-\alpha)\sqrt{\gamma+\alpha^2-1}},D_0+\frac{\delta}{R(1-\alpha)\sqrt{\gamma+\alpha^2-1}}\right)$
to a~neighborhood $\left(D_1-\frac{\delta}{R(1+\alpha)\sqrt{\gamma+\alpha^2-1}},D_1+\frac{\delta}{R(1+\alpha)\sqrt{\gamma+\alpha^2-1}}\right)$,
the minimum of $V(\cos\theta,D,\lambda)$ moves from a~neighborhood of $\theta=0$ to a~neighborhood of $\theta=\pi$.

As explained before, it is suf\/f\/icient to study the position of $z_{\min}$ for the function
$2I_3R^2\gamma^2(\gamma+\alpha^2-1)\tilde{V}(z,D(a,b),\lambda(a,b))=:f(z)$ with
\begin{gather*}
f(z)=\frac{(az+b)^2}{1-z^2}-\beta z.
\end{gather*}
The derivative is $f'(z)=\frac{p(z)}{(1-z^2)^2}$, where $p(z)=2(az+b)(bz+a)-\beta(1-z^2)^2$.

Our question about the position of $z_{\min}$ can be reduced to the question: for which values of the
parameters $a,b,\beta$ does the polynomial $p(z)$ have one zero in one of the subintervals $[1-\epsilon,1]$
and $[-1,-1+\epsilon]$ of the interval $[-1,1]$?

We see that $p(-1)=-2(a-b)^2$ and $p(1)=2(a+b)^2$, which conf\/irms that $f'(z)$ has a~zero and changes
sign in $[-1,1]$ (remember that $a\neq\pm b$ if $D\neq D_1$ and $D\neq D_0$).
\begin{proposition}\label{local_min}
Assume that
$\lambda>\lambda_{\text{\rm thres}}=\frac{\sqrt{mgR^3I_3\alpha}(1+\alpha)^2}{\sqrt{1+\alpha-\gamma}}$. %\text~31
\begin{itemize}\itemsep=0pt
\item[$i)$] For any $($small$)$ $\epsilon>0$ there is $\delta_{-}(\epsilon,\lambda)>0$ being a~positive solution
of the equation
\begin{gather*}
\alpha(1-\alpha)\gamma\lambda\delta+\alpha\delta^2+\frac{\delta^2}{\epsilon^2}
(1-\epsilon)(1+\alpha)^2=\alpha^2\gamma^2\lambda^2-\frac{\beta}{2}(2-\epsilon)^2(1+\alpha)^2,
\end{gather*}
such that for every positive $\delta<\delta_{-}(\epsilon,\lambda)$ the potential $V(z,D,\lambda)$ has
a~minimum $z_{\min}$ in the interval $[-1,-1+\epsilon]$ for
$D=D_1\pm\frac{\delta}{R(1+\alpha)\sqrt{\gamma+\alpha^2-1}}$.

\item[$ii)$] For any $($small$)$ $\epsilon>0$ there is
$\delta_{+}(\epsilon,\lambda)=\min\{\delta_1,\delta_2\}>0$, where ${\delta_1=\gamma(1+\alpha)\lambda}$ and
$\delta_2$ is the positive solution of the equation
\begin{gather*}
\frac{\delta^2}{\epsilon^2}
(1-\epsilon)(1-\alpha)^2+\alpha\gamma\lambda\delta(1+\alpha)-\alpha\delta^2=\alpha^2\gamma^2\lambda^2+\frac{\beta}{2}(2-\epsilon)^2(1-\alpha)^2,
\end{gather*}
such that for every positive $\delta<\delta_{+}(\epsilon,\lambda)$ the potential $V(z,D,\lambda)$ has
a~minimum $z_{\min}$ in the interval $[1-\epsilon,1]$ for
$D=D_0\mp\frac{\delta}{R(1-\alpha)\sqrt{\gamma+\alpha^2-1}}$.
\end{itemize}
\end{proposition}
\begin{proof}
We know that $V(z,D,\lambda)$ has one minimum in $[-1,1]$ and $p(z)$, the numerator of the derivative
$f'(z)$, has one zero in $[-1,1]$.
The question is to formulate conditions for the para\-me\-ter~$D$ so that this zero is located in the interval
$[-1,-1+\epsilon]$ or $[1-\epsilon,1]$.

 $i)$ We have to show that for any (small) $\epsilon>0$ and $\lambda>\lambda_{\text{thres}}$ there
is ${\delta_{-}(\epsilon,\lambda)>0}$ such that for any $\delta<\delta_{-}(\epsilon,\lambda)$ the potential
$V\Big(z,D=D_1 \pm\frac{\delta}{R(1+\alpha)\sqrt{\gamma+\alpha^2-1}},\lambda\Big)$ has a~minimum in the interval $[-1,-1+\epsilon]$.
To f\/ind this $\delta_{-}(\epsilon,\lambda)$ we have to study the signs of $p(-1)=-2(a-b)^2<0$ and of
\begin{gather*}
p(-1+\epsilon)=p(-1)+2(a-b)^2\epsilon+2(ab-2\beta)\epsilon^2+4\beta\epsilon^3-\beta\epsilon^4.
\end{gather*}
To have minimum $z_{\min}$ in $[-1,-1+\epsilon]$ we need to have $p(-1+\epsilon)>0$.
If $D=D_1\pm\frac{\delta}{R(1+\alpha)\sqrt{\gamma+\alpha^2-1}}$ then $(a-b)^2=\delta^2$,
$ab=\frac{1}{(1+\alpha)^2}(\alpha\gamma\lambda-\delta)(\alpha\gamma\lambda+\alpha\delta)$ and
\begin{gather*}
p(-1+\epsilon)=-2\delta^2+2\delta^2\epsilon+2ab\epsilon^2-\beta\epsilon^2(2-\epsilon)^2
\\
\phantom{p(-1+\epsilon)}
=-2\delta^2(1-\epsilon)+\frac{2\epsilon^2}{(1+\alpha)^2}
\left(\alpha^2\gamma^2\lambda^2-\alpha\gamma(1-\alpha)\delta\lambda-\alpha\delta^2\right)-\beta\epsilon^2(2-\epsilon)^2.
\end{gather*}
For positivity of $p(-1+\epsilon)$ we must have
\begin{gather}
\label{delta-1}
\alpha(1-\alpha)\gamma\lambda\delta+\alpha\delta^2+\frac{\delta^2}{\epsilon^2}
(1-\epsilon)(1+\alpha)^2<\alpha^2\gamma^2\lambda^2-\frac{\beta}{2}(2-\epsilon)^2(1+\alpha)^2.
\end{gather}
The r.h.s.
of~\eqref{delta-1} is positive since
\begin{gather*}
(\alpha\gamma\lambda)^2-\frac{\beta}{2}
(2-\epsilon)^2(1+\alpha)^2=(\alpha\gamma\lambda)^2\left[1-\frac{2\beta}{\lambda^2}\frac{(1+\alpha)^2}
{\alpha^2\gamma^2}\frac{(2-\epsilon)^2}{4}\right],
\end{gather*}
and for $\lambda=C\lambda_{\text{thres}}$, $C>1$, the term in the square brackets can be shown to be positive.
We see this when we write
\begin{gather}
\label{gamma_alpha_estimate}
\frac{2\beta}{\lambda^2}\frac{(1+\alpha)^2}{\alpha^2\gamma^2}\frac{(2-\epsilon)^2}{4}
=4\left[\frac{(\gamma+\alpha^2-1)(1+\alpha-\gamma)}{\alpha^2(1+\alpha)^2}\right]\frac{1}{C^2}
\frac{(2-\epsilon)^2}{4}.
\end{gather}
Clearly $\frac{1}{C^2}\frac{(2-\epsilon)^2}{4}<1$ and the function
$\rho(\alpha,\gamma)=\frac{(\gamma+\alpha^2-1)(1+\alpha-\gamma)}{\alpha^2(1+\alpha)^2}<\frac{1}{4}$ for
${1-\alpha^2<\gamma<1}$, $0<\alpha<1$.
So the r.h.s.
of~\eqref{delta-1} is indeed positive.
The l.h.s.
of~\eqref{delta-1} has the form
$\delta[\alpha(1-\alpha)\gamma\lambda+\alpha\delta+\frac{\delta}{\epsilon^2}(1-\epsilon)(1+\alpha)^2]$ with
all terms positive.
If we def\/ine $\delta_{-}(\epsilon,\lambda)$ as the positive solution of
\begin{gather}
\label{delta-1-1}
\alpha(1-\alpha)\gamma\lambda\delta+\alpha\delta^2+\frac{\delta^2}{\epsilon^2}
(1-\epsilon)(1+\alpha)^2=\alpha^2\gamma^2\lambda^2-\frac{\beta}{2}(2-\epsilon)^2(1+\alpha)^2,
\end{gather}
then for all positive $\delta<\delta_{-}(\epsilon,\lambda)$ we have $p(-1+\epsilon)>0$ and the minimum
$z_{\min}\in[-1,-1+\epsilon]$.

  $ii)$ We have to show that for any (small) $\epsilon>0$ and $\lambda>\lambda_{\text{thres}}$
there is $\delta_{+}(\epsilon,\lambda)>0$ such that for any positive $\delta<\delta_{+}(\epsilon,\lambda)$
the potential ${V\Big(z,D=D_{0}\mp\frac{\delta}{R(1-\alpha)\sqrt{\gamma+\alpha^2-1}}\Big)}$ has
a~minimum $z_{\min}$ in the interval $[1-\epsilon,1]$.
To f\/ind this $\delta_{+}(\epsilon,\lambda)$ we need to study the signs of $p(1)=2(a+b)^2>0$ and of
\begin{gather*}
p(1-\epsilon)=p(1)-2(a+b)^2\epsilon+2(ab-2\beta)\epsilon^2+\beta\epsilon^3(4-\epsilon).
\end{gather*}
To have minimum $z_{\min}$ in $[1-\epsilon,1]$ we need to show how to f\/ind $\delta$ so that
${p(1-\epsilon)<0}$.
If $D=D_0\mp\frac{\delta}{R(1-\alpha)\sqrt{\gamma+\alpha^2-1}}$ then $(a+b)^2=\delta^2$,
${ab=-\frac{1}{(1-\alpha)^2}(\alpha\gamma\lambda-\delta)(\alpha\gamma\lambda-\alpha\delta)}$ and
\begin{gather*}
p(1-\epsilon)=2\delta^2(1-\epsilon)-\frac{2\epsilon^2}{(1-\alpha)^2}
((\alpha\gamma\lambda-\delta)(\alpha\gamma\lambda-\alpha\delta)+2\beta(1-\alpha)^2)+\beta\epsilon^3(4-\epsilon).
\end{gather*}
For negativity of $p(1-\epsilon)$ we must have
\begin{gather}
\label{delta+1}
\frac{\delta^2}{\epsilon^2}
(1-\epsilon)(1-\alpha)^2+\alpha\gamma\lambda\delta(1+\alpha)-\alpha\delta^2<\alpha^2\gamma^2\lambda^2+\frac{\beta}{2}(2-\epsilon)^2(1-\alpha)^2.
\end{gather}
The r.h.s.\
of~\eqref{delta+1} is obviously positive.
The l.h.s.\
of~\eqref{delta+1} has the form
\begin{gather*}
{\delta\left[\frac{\delta}{\epsilon^2}(1-\epsilon)(1-\alpha)^2+\alpha((1+\alpha)\gamma\lambda-\alpha\delta)\right]},
\end{gather*}
which is certainly positive when ${\lambda>\frac{\delta}{\gamma(1+\alpha)}}$ since $\delta>0$
can always be chosen so small that this is true, by the chosen convention $\lambda>0$.
If we def\/ine $\delta_{+}(\epsilon,\lambda)=\min\{\delta_1,\delta_2\}$, where
$\delta_1=\gamma(1+\alpha)\lambda$ and $\delta_2$ is the positive solution of
\begin{gather*}
\frac{\delta^2}{\epsilon^2}
(1-\epsilon)(1-\alpha)^2+\alpha\gamma\lambda\delta(1+\alpha)-\alpha\delta^2=\alpha^2\gamma^2\lambda^2+\frac{\beta}{2}(2-\epsilon)^2(1-\alpha)^2,
\end{gather*}
then for all positive $\delta<\delta_{+}(\epsilon,\lambda)$ we have $p(1-\epsilon)<0$ and the minimum
${z_{\min}\in[1-\epsilon,1]}$.
\end{proof}
\begin{example}
\label{main_example}
We exemplify the above proposition for some physical values for the parameters.
We let $m=0.02$ kg, $R=0.02$ m and $\alpha=0.3$.
We shall assume $I_3=\frac{2}{5}mR^2$, and we get a~rational potential~\eqref{rational_pot} when
$I_1=\frac{I_3(I_3+mR^2(1-\alpha^2))}{I_3+mR^2}=\frac{131}{350}mR^2$, which means that
$\gamma=\frac{131}{140}\in(1-\alpha^2,1)=(0.91,1)$.
Jellett's integral $\lambda$ will for simplicity be assumed to be two times the threshold value
$\lambda_{\text{thres}}$, so $\lambda=6.88\cdot 10^{-6}$ $\text{kg}\cdot\text{m}^3\cdot\text{s}^{-1}$.

Suppose f\/irst that $\epsilon=0.1$ and consider the condition for a~minimum $z_{\min}$ to stay in the
interval $(-1,-1+\epsilon)$.
We have $D_1=-6\cdot 10^{-4}$ and by solving~\eqref{delta-1-1} we get $\delta_{-}=1.48\cdot 10^{-7}$, so if
$D=D_1+\delta$, with $\delta<\frac{\delta_{-}}{R(1+\alpha)\sqrt{\gamma+\alpha^2-1}}\approx 3.5\cdot
10^{-5}$, the potential has its minimum in the interval.

If $\epsilon$ is reduced to $0.01$, the bound for $\delta$ is tightened one order of magnitude as well.
\end{example}

\section[Oscillation of $\theta(t)$ within the deforming rational potential ${V(\cos\theta,D(t),\lambda)}$]{Oscillation
of $\boldsymbol{\theta(t)}$ within the deforming rational\\ potential $\boldsymbol{V(\cos\theta,D(t),\lambda)}$} %\\

As a~toy TT inverts, we can see that the symmetry axis performs oscillations, or equivalently we say that
it nutates.
This is also apparent in simulations of the equations of motion~\cite{Coh,Or,Ued} where graphs of the
evolution of the inclination angle $\theta(t)$ show that it rises in an oscillating manner from an angle
close to $\theta=0$ to an angle close to $\theta=\pi$.
Here we say that a~solution $\theta(t)$ is oscillatory on a~time interval $[0,T]$ when $\dot{\theta}(t)$
changes sign a~number of times in this interval.

We consider solutions of METT with $D(t)$, $\tilde{E}(t)$ describing inverting solutions of TT equations,
under assumption that $D(t)$, $\tilde{E}(t)$ are slowly varying regular functions moving from a~small
neighborhood of $(D_0,\tilde{E}_0)$ to a~small neighborhood of $(D_1,\tilde{E}_1)$.
They are regular since $D(t)=D(\theta(t),\omega_3(t))$ and
$\tilde{E}(t)=\tilde{E}(\theta(t),\dot{\theta}(t),\dot{\varphi}(t),\omega_3(t))$.
We further have assumed that $(D(t),\tilde{E}(t))\stackrel{t\to\infty}{\longrightarrow}(D_1,\tilde{E}_1)$
and $(D(t),\tilde{E}(t))$ ${\stackrel{t\to-\infty}{\longrightarrow}(D_0,\tilde{E}_0)}$.

In the limiting case of constant $D$ and $\tilde{E}$ when the METT describes purely rolling TT, the
oscillating behaviour of $\theta(t)$ follows from the dynamical system representation of the second order
equation for $\theta(t)$.
For the rolling TT the energy
\begin{gather}
\label{METTsection5}
\tilde{E}=g(\cos\theta)\dot{\theta}^2+V(\cos\theta,D,\lambda)
\end{gather}
is an integral of motion for the $\theta$-equation obtained by dif\/ferentiating~\eqref{METTsection5}:
\begin{gather*}
\ddot{\theta}=\frac{\sin\theta}{2g(\cos\theta)}\left(g_{z}'(\cos\theta)\dot{\theta}^2+V_{z}'(\cos\theta,D,\lambda)\right),
\end{gather*}
where $g(z=\cos\theta)=\frac{I_3(\alpha^2+(1-\gamma)^2-2\alpha(1-\gamma)z)}{2(\gamma+\alpha^2-1)}$.
This means that $\tilde{E}=g(\cos\theta)y^2+V(\cos\theta,D,\lambda)$ is an integral of motion for the
dynamical system
\begin{gather}
\dot{\theta}=y,
\nonumber
\\
\dot{y}=\frac{\sin\theta}{2g(\cos\theta)}\left(g_{z}'(\cos\theta)y^2+V_{z}'(\cos\theta,D,\lambda)\right).
\label{dyn_sys_theta}
\end{gather}
Trajectories of system~\eqref{dyn_sys_theta} are lines of constant value of energy $\tilde{E}$ and they are
closed curves.
The closed trajectories describe periodic solutions~\cite{Chap2} with period $T$ def\/ined by the integral
\begin{gather*}%\label{nutation_period_ish}
T(\tilde{E})=2\int_{\theta_1}^{\theta_2}\frac{\sqrt{g(\cos\theta)}d\theta}{\sqrt{\tilde{E}-V(\cos\theta,D,\lambda)}},
\end{gather*}
where the turning latitudes $\theta_1<\theta_2$ are def\/ined by $\tilde{E}=V(\cos\theta_{1,2},D,\lambda)$.

The modif\/ied energy $\tilde{E}(t)$ is bounded.
That entails boundedness of
$\boldsymbol{\omega}(t)=-\dot{\varphi}\sin\theta\mathbf{\hat{1}}+\dot{\theta}\mathbf{\hat{2}}+(\dot{\psi}+\dot{\varphi}\cos\theta)\mathbf{\hat{3}}$
and thus $|\dot{\theta}(t)|<B$ for some positive $B$.
Since the potential ${V(z,D(t),\lambda)\stackrel{\theta\to 0,\pi}{\longrightarrow}\infty}$ the curve
$(\theta(t),\dot{\theta}(t))$ of each inverting solution is conf\/ined to the open rectangle
$(\theta,\dot{\theta})\in(0,\pi)\times(-B,B)$.

The picture that emerges is that, for slowly varying $\tilde{E}(t)$, $D(t)$, the inverting trajectories of
METT stay (locally) close to the trajectories of system~\eqref{dyn_sys_theta} and are traversed with almost
the same velocity.
The time of passing $T=t_2-t_1$ between two turning angles given by
$V(\cos\theta_1,D(t_1),\lambda)=\tilde{E}(t_1)$ and $V(\cos\theta_2,D(t_2),\lambda)=\tilde{E}(t_2)$ is
close to the half-period of the non-deforming potential.

Thus initially the trajectory moves around $(\theta_{\min},\dot{\theta}=0)$, with $\theta_{\min}$ close to~0.
As $D(t)\to D_1$ the minimum $\theta_{\min}$ moves toward $\theta=\pi$ (see Proposition~\ref{local_min})
and for suf\/f\/iciently slowly va\-rying~$D(t)$ the trajectory goes several times around
($\theta_{\min}$, $\dot{\theta}=0$) and drifts toward the point ($\theta=\pi$, $\dot{\theta}=0$).
The trajectory describes a~selfcrossing and contracting spiral with the center moving toward $(\pi,0)$ in
the $(\theta,\dot{\theta})$-plane.
\begin{figure}[t]
\centering
\includegraphics[scale=0.5]{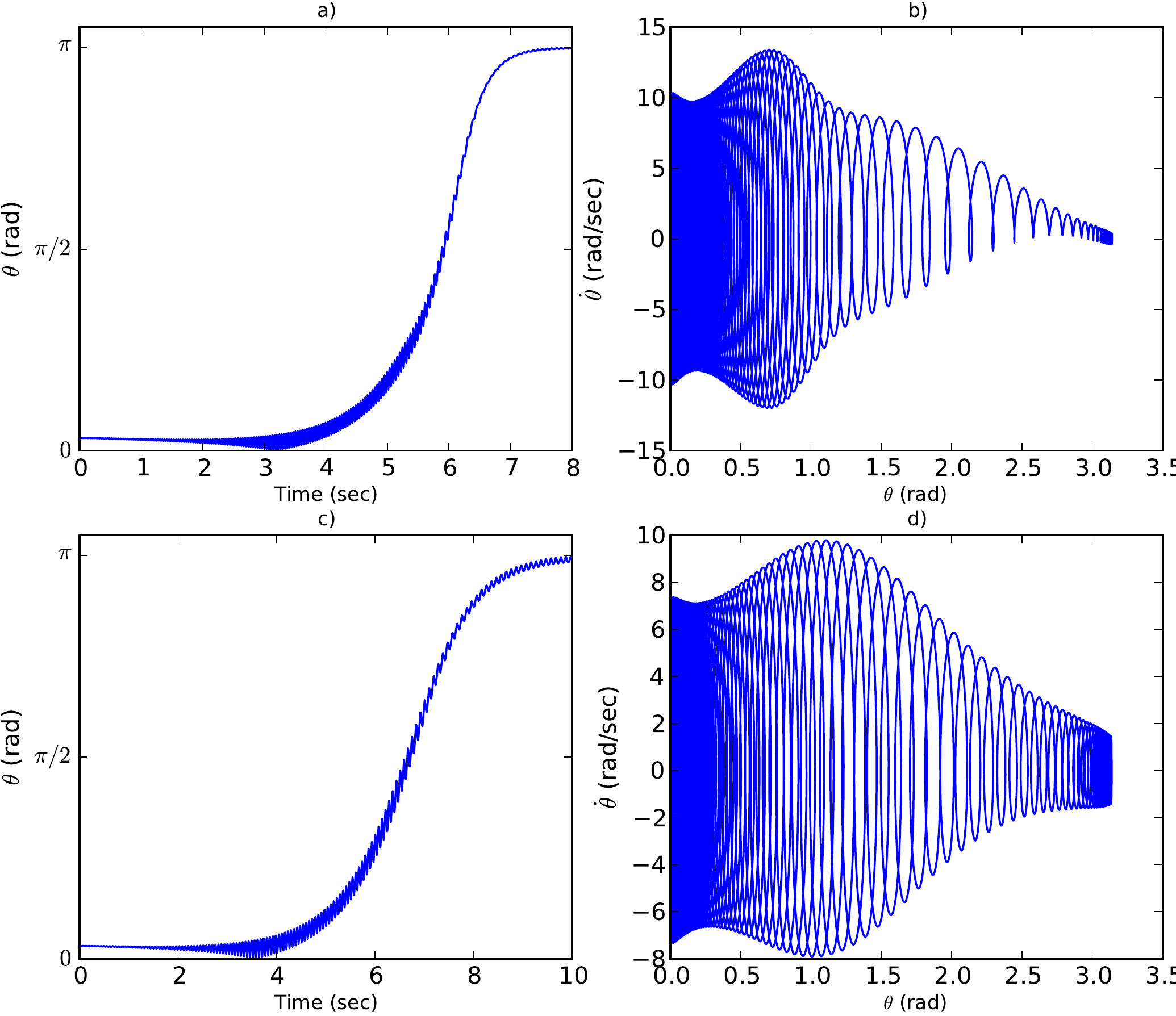}

\caption{Plots of $(t,\theta(t))$ (left) and
$(\theta(t),\dot{\theta}(t))$ (right) obtained by integrating equations~\eqref{ddth3}--\eqref{nu_y3}
and~\eqref{equ_for_g_n} for two sets of parameters and initial values.
Plots $a$ and $b$ correspond to parameter values for rational potential in Example~\ref{main_example} and
$\mu=0.3$.
Plots $c$ and $d$ correspond to parameter values provided in~\cite{Coh}, corresponding to an algebraic
potential.
The equations are integrated using the Python~2.7 open source library SciPy~\cite{Scipy}.}\label{quadplot}
\end{figure}
The solution curve for $\theta(t)$ is a~superposition of a~small amplitude oscillating component with
a~logistic type curve going from $\theta=0$ to $\theta=\pi$.

In Fig.~\ref{quadplot} we see this.
The plots $a$ and $b$ present a~graph of $\theta(t)$ and a~curve $(\theta(t),\dot{\theta}(t))$ obtained by
solving equations of motion~\eqref{ddth3}--\eqref{nu_y3} for the rolling and gliding TT, with $g_n$ given
by equation~\eqref{equ_for_g_n} and with values of parameters taken from Example~\ref{main_example}.
So the ef\/fective potential~$V$ given by~\eqref{rational_pot} is rational.
We have taken~$\mu=0.3$ and we have integrated the equations numerically for initial values
$\theta(0)=0.1\text{~rad}$, $\dot{\theta}(0)=\dot{\varphi}(0)=0$, $\nu_x(0)=\nu_y(0)=0$ and
$\omega_3(0)=155.0\text{~rad/s}$ (corresponding to the value of $\lambda=2\lambda_{\text{thres}}$).

Plot $a$ shows oscillations of $\theta(t)$ as it rises to $\pi$.
The time of inversion from the moment when~$\theta(t)$ starts to rise, at approximately 3 seconds, is about
4--5 seconds.
Larger values of $\mu$ give shorter time of inversion.
Plot $b$ shows the trajectory in the $(\theta,\dot{\theta})$-plane.

In Fig.~\ref{quadplot} $c$ and $d$ the same type of numerical integration is performed,
for parameters giving an algebraic (not rational) ef\/fective potential.
The parameters are the same as used in Cohen~\cite{Coh} and they are close to values taken in
Example~\ref{main_example}.
Here $m=0.015\text{ kg}$, $R=0.025\text{~m}$, $\alpha=0.2$, $I_1=I_3=\frac{2}{5}mR^2$ and $\mu=0.3$ with
initial values $\theta(0)=0.1\text{~rad}$, $\dot{\theta}(0)=\dot{\varphi}(0)=0$, $\nu_x(0)=\nu_y(0)=0$ and
$\omega_3(0)=100\text{ rad/s}$.
We may note that Cohen uses a~Coulomb-type gliding friction force, so the equations solved in his article
are slightly dif\/ferent than equations~\eqref{ddth3}--\eqref{nu_y3}.

In the next section we shall estimate the maximum value of the period of oscillations for all values of the
modif\/ied energy $\tilde{E}$ and $D$.

This will allow us to formulate a~condition for the time of inversion needed for oscillatory behaviour of
$\theta(t)$.
Basically $T_{\text{inv}}$ has to be an order of magnitude larger than $T_{\text{upp}}$~-- the maximal
period of oscillations within potential $V(\cos\theta,D(t),\lambda)$. %\text

\subsection{Estimates for the period of oscillation}

A direct way of estimating the period of nutation is to study the explicit integral def\/ining the period
and expanding it w.r.t.\
a small parameter $\epsilon=\frac{2\beta}{b^2}$, with $b$ and $\beta$ given by~\eqref{b_def}
and~\eqref{beta_def}, respectively.
This is the technique used in~\cite{Nisse3}, but in that paper the analysis is based on the assumption that
$\omega_3$ is large.
Here we use the more general assumption that $\lambda$ is only above the threshold value, so that
$\lambda=C\lambda_{\text{thres}}$ with $C>1$.
To simplify estimates we also assume that we consider curves $(D(t),\tilde{E}(t))$ such that $D_1<D<D_0$.
Then by using that $b=\alpha\lambda+\alpha RD\sqrt{\gamma+\alpha^2-1}$ satisf\/ies
$\frac{\alpha\gamma\lambda}{1+\alpha}<b<\frac{\alpha\gamma\lambda}{1-\alpha}$ for $D_1<D<D_0$, we have the
estimate
\begin{gather}
\label{epsilon_estimate}
\epsilon=\frac{2\beta}{b^2}<\frac{2\beta(1+\alpha)^2}{\lambda^2\alpha^2\gamma^2}=\frac{4}{C^2}
\cdot\left[\frac{(\gamma+\alpha^2-1)(1+\alpha-\gamma)}{\alpha^2(1+\alpha)^2}\right]<\frac{4}{C^2}
\cdot\frac{1}{4}<1,
\end{gather}
for $1-\alpha^2<\gamma<1$, $0<\alpha<1$, similar to the estimate for equation~\eqref{gamma_alpha_estimate}.

In the following we shall determine the dependence of the period $T(\epsilon)$ on the small parameter
$\epsilon=\frac{2\beta}{b^2}$ and we shall f\/ind an estimate for the maximal value of $T(\epsilon)$ that
is valid for all $D\in(D_1,D_0)$ and all values $\epsilon<1$, which means all values of the Jellett
integral that are above the threshold value.
As known from the asymptotic analysis these initial values of $\lambda$ lead to inversion of the TT.

The period of oscillations for the potential $V(z,D,\lambda)$ is given by the integral
\begin{gather}
\label{Integral_for_period}
T=2\int_{z_1}^{z_2}\frac{\sqrt{g(z)}dz}{\sqrt{(1-z^2)(\tilde{E}-V(z,D,\lambda))}},
\end{gather}
where $z_1<z_2\in(-1,1)$ are two turning points def\/ined by $\tilde{E}=V(z_{1,2},D,\lambda)$.
This equation always has two solutions for $\tilde{E}>V(z_{\min},D,\lambda)$ since the potential is convex
and %${}$
\begin{gather*}
\lim_{z\to\pm 1}V(z,D,\lambda)=\infty.
\end{gather*}

In terms of the the parameters $a$ and $b$, the potential~\eqref{rational_pot} reads
\begin{gather*}
V(z,D,\lambda)=mgR(1-\alpha z)+\frac{(az+b)^2}{2I_3R^2\gamma^2(\gamma+\alpha^2-1)(1-z^2)}
+\frac{\alpha^2a^2-(1-\gamma)b^2}{2R^2\alpha^2I_3\gamma^2(\gamma+\alpha^2-1)},
\end{gather*}
and (by def\/inition) the left turning point $z_1$ is given implicitly by the equation
\begin{gather*}%\label{Energy_at_z1}
\tilde{E}=mgR(1-\alpha z_1)+\frac{(az_1+b)^2}{2I_3R^2\gamma^2(\gamma+\alpha^2-1)(1-z_{1}^{2})}
+\frac{\alpha^2a^2-(1-\gamma)b^2}{2R^2\alpha^2I_3\gamma^2(\gamma+\alpha^2-1)}.
\end{gather*}
In the following we shall parametrise (similarly as in~\cite{Nisse3}) the remaining roots of
$(1-z^2)(\tilde{E}-V(z,D,\lambda))$ by $z_1$.
Due to this they become solutions of a~quadratic equation.
The function in the denominator of~\eqref{Integral_for_period} we write as
\begin{gather}
\big(1-z^2\big)(\tilde{E}-V(z,D,\lambda))=\frac{1}{\big(1-z_{1}^2\big)}\biggl(-mgR\alpha(z_1-z)(1-z_{1}^2)(1-z^2)
\nonumber
\\
\quad\qquad
{}+\frac{(az_1+b)^2(1-z^2)-(az+b)^2(1-z_{1}^2)}{2I_3R^2\gamma^2(\gamma+\alpha^2-1)}\biggr)
\nonumber
\\
\qquad
=\frac{\beta(z_1-z)}{2I_3R^2\gamma^2(\gamma+\alpha^2-1)}\left(z^2+\frac{a^2+b^2+2abz_1}{\beta(1-z_{1}^2)}
z-1+\frac{(a^2+b^2)z_1+2ab}{\beta(1-z_{1}^2)}\right).
\label{equz1z2z3}
\end{gather}
Notice that $z=\pm1$ are {\it not} roots because $V(z,D,\lambda)$ has singularities at the points $z=\pm 1$.
The quadratic polynomial in the parentheses of~\eqref{equz1z2z3} determines roots $z_2$, $z_3$, for any
given turning point $z_1$.
They are thus solutions of the quadratic equation
\begin{gather}
\label{z2z3_equ}
z^2+\left(\frac{a^2+b^2+2abz_1}{\beta(1-z_{1}^2)}\right)z-1+\frac{(a^2+b^2)z_1+2ab}{\beta(1-z_{1}^2)}=0.
\end{gather}
By $z_3$ we denote the root satisfying $z_3<-1$ and by $z_2$ $({>}z_1)$ the right turning point.
The polynomial $(1-z^2)(E-V(z,D,\lambda))$ is then factorized:
\begin{gather*}
\big(1-z^2\big)(E-V(z,D,\lambda))=\frac{\beta(z_1-z)(z_2-z)(z_3-z)}{2I_3R^2\gamma^2(\gamma+\alpha^2-1)}
=mgR\alpha(z_1-z)(z_2-z)(z_3-z).
\end{gather*}
We substitute this into the integral for the period $T$~\eqref{Integral_for_period}.
By the mean value theorem for integrals there exist a~$z^*\in[z_1,z_2]$ such that we have
\begin{gather*}
T=\frac{2\sqrt{g(z^{*})}}{\sqrt{mgR\alpha}}\int_{z_1}^{z_2}\frac{dz}{\sqrt{(z_1-z)(z_2-z)(z_3-z)}},
\end{gather*}
where $g(z^{*})=\frac{I_3(\alpha^2+(1-\gamma)^2-2\alpha(1-\gamma)z^{*})}{2(\gamma+\alpha^2-1)}$.

The integral can be transformed to a~standard complete elliptic integral of the f\/irst kind through the
change of variables $z=z_2+(z_1-z_2)s^2$~\cite{Bate}:
\begin{gather}
\label{period_integral}
T=\frac{4\sqrt{g(z^{*})}}{\sqrt{mgR\alpha}}\frac{1}{\sqrt{z_2-z_3}}\int_{0}^{1}\frac{ds}
{\sqrt{(1-s^2)(1-k^2s^2)}},
\end{gather}
where $k^2=\frac{z_2-z_1}{z_2-z_3}<1$ is a~positive parameter.
This integral has the standard expansion
\begin{gather}
\label{complete_elliptic}
K\big(k^2\big)=\int_{0}^{1}\!\!\frac{ds}{\sqrt{(1-s^2)(1-k^2s^2)}}=\frac{\pi}{2}\sum_{n=0}^{\infty}
\left(\frac{(2n-1)!!}{(2n)!!}\right)^2\! k^{2n}=\frac{\pi}{2}\left(1+\frac{k^2}{4}+O\big(k^4\big)\right).\!\!\!
\end{gather}
By using the roots of equation~\eqref{z2z3_equ} and by expanding $k^2=\frac{z_2-z_1}{z_2-z_3}$ w.r.t.\
 $\epsilon=\frac{2\beta}{b^2}$, we can show that $k^2=O(\epsilon)$ if $\epsilon$ is small.
Indeed, after solving equation~\eqref{z2z3_equ} we get that the quantity $k^2$ can be written:
\begin{gather}
k^2=\frac{z_2-z_1}{z_2-z_3}=\frac{1}{2}-\frac{1}{2}\left(z_1+\frac{a^2+b^2+2abz_1}{2\beta(1-z_{1}^2)}
\right)\nonumber\\
\phantom{k^2=}{}\times\left(\frac{(a^2+b^2+2abz_1)^2}{4\beta^2(1-z_{1}^2)^2}+1-\frac{(a^2+b^2)z_1+2ab}{\beta(1-z_{1}^2)}
\right)^{-\frac{1}{2}}
\nonumber
\\
\phantom{k^2}
=\frac{1}{2}-\frac{1}{2}\left(z_1\frac{2\beta\big(1-z_{1}^2\big)}{a^2+b^2+2abz_1}
+1\right)
\nonumber\\
\phantom{k^2=}{}\times
\left(1+\frac{4\beta^2(1-z_{1}^2)^2}{(a^2+b^2+2abz_1)^2}
-\frac{4\beta(1-z_{1}^2)((a^2+b^2)z_1+2ab)}{(a^2+b^2+2abz_1)^2}\right)^{-\frac{1}{2}}.
\label{k2_3}
\end{gather}
With the parameters $\epsilon=\frac{2\beta}{b^2}$ and $w=\frac{a}{b}$, the expansion of $k^2$ is:
\begin{gather}
k^2=\frac{1}{2}-\frac{1}{2}\left(\frac{z_1(1-z_{1}^2)}{w^2+1+2wz_1}
\epsilon+1\right)\nonumber\\
\phantom{k^2=}{}\times
\left(1+\frac{(1-z_{1}^2)^2}{(w^2+1+2wz_1)^2}\epsilon^2-\frac{2(1-z_{1}
^2)((w^2+1)z_1+2w)}{(w^2+1+2wz_1)^2}\epsilon\right)^{-\frac{1}{2}}
\nonumber
\\
\phantom{k^2}
=-\frac{(1+z_1w)(w+z_1)\big(1-z_{1}^2\big)}{(1+w^2+2wz_1)^2}\epsilon+O\big(\epsilon^2\big).
\label{k2_est3}
\end{gather}
It should be noted that the minus sign here is misleading.
When~\eqref{k2_est3} is expressed by parame\-ters~$D$,~$\lambda$, one sees that the factor at $\epsilon$ is
positive.

We consider the parameter $k^2$ for $(z_1,w)\in[-1,1]\times[-1+\delta,1-\delta]$ with certain small
$\delta$.
The values $w=\frac{a}{b}=\mp1$ correspond (as Fig.~\ref{ab_diagram} shows) to the upright and inverted
spinning solutions, which are asymptotic solutions of TT equations and are never attained during the
inversion.

The inverting TT starts with $w_0=-1+\delta_0$ and angle $\theta_0$ close to 0 and moves to the value
$w_1=1-\delta_1$ where the angle $\theta_1$ is close to $\pi$.
So we need to estimate the function $h_{1}(z_1,w)=-\frac{(1+z_1 w)(w+z_1)(1-z_{1}^2)}{(1+w^2+2wz_1)^2}$,
being a~factor in front of $\epsilon$ in~\eqref{k2_est3}, on a~rectangle
$[-1,1]\times[-1+\delta,1-\delta]=:R_{\delta}$ with $\delta<\min\{\delta_0,\delta_1\}$.
The function $h_1(z_1,w)$ has on $[-1,1]\times(-1,1)$ two critical points $(\frac{1}{\sqrt{3}},0)$,
$(-\frac{1}{\sqrt{3}},0)$, which gives a~local minimum $h_1(\frac{1}{\sqrt{3}},0)=-\frac{2}{3\sqrt{3}}$ and
a~local maximum $h_{1}(-\frac{1}{\sqrt{3}},0)=\frac{2}{3\sqrt{3}}$.
In the whole region $R_{\delta}$ the function $h(z_1,w)$ satisf\/ies the inequality
${-\frac{2}{3\sqrt{3}}\leq h_{1}(z_1,w)\leq\frac{2}{3\sqrt{3}}}$ and, therefore, it is bounded
$|h_{1}(z_1,w)|\leq\frac{2}{3\sqrt{3}}$.
This clarif\/ies how $k^2$ depends on the small parameter $\epsilon$.

To estimate $\frac{1}{\sqrt{z_2-z_3}}$ in~\eqref{period_integral}, we write
\begin{gather}
\frac{1}{\sqrt{z_2-z_3}}=\frac{1}{\sqrt{2}}\left(\frac{\big(a^2+b^2+2abz_1\big)^2}{4\beta^2(1-z_{1}^2)^2}
+1-\frac{((a^2+b^2)z_1+2ab)}{\beta(1-z_{1}^2)}\right)^{-1/4}\!\!=\frac{\sqrt{\beta\big(1-z_{1}^2\big)}}{b\sqrt{1+w^2+2wz_1}}
\nonumber
\\
\phantom{\frac{1}{\sqrt{z_2-z_3}}=}{}\times
\left(1+\frac{(1-z_{1}^2)^2}{(1+w^2+2wz_1)^2}
\epsilon^2-\frac{2(1-z_{1}^2)((1+w^2)z_1+2w)}{(1+w^2+2wz_1)^2}\epsilon\right)^{-1/4}
\nonumber
\\
\phantom{\frac{1}{\sqrt{z_2-z_3}}}
=\sqrt{\frac{\epsilon}{2}}\frac{\sqrt{\big(1-z_{1}^2\big)}}{\sqrt{1+w^2+2wz_1}}\left(1+\frac{(1-z_{1}
^2)((1+w^2)z_1+2w)}{2(1+w^2+2wz_1)^2}\epsilon+O\big(\epsilon^2\big)\right).
\label{sqrtz2z3}
\end{gather}
The function $h_{2}(z_1,w)=\frac{(1-z_{1}^2)((1+w^2)z_1+2w)}{2(1+w^2+2wz_1)^2}$ satif\/ies the inequalities
${-1<h_2(z_1,w)<1}$ on the rectangle $(z_1,w)\in R_{\delta}$, so it is bounded: $|h_2(z_1,w)|\leq1$.

Thus the nutational period behaves as
\begin{gather*}
T(\epsilon)=\frac{4\sqrt{g(z^*)}}{\sqrt{mgR\alpha}}\frac{1}{\sqrt{z_2-z_3}}\int_{0}^{1}\frac{ds}
{\sqrt{(1-s^2)(1-k^2s^2)}}
\\
\phantom{T(\epsilon)}
=\frac{4\sqrt{g(z^*)}}{\sqrt{mgR\alpha}}\sqrt{\frac{\epsilon}{2}}\frac{\sqrt{1-z_{1}^2}}
{\sqrt{1+w^2+2wz_1}}\left[1+h_{2}(z_1,w)\epsilon+O\big(\epsilon^2\big)\right]
\\
\phantom{T(\epsilon)=}{}\times
\frac{\pi}{2}\left[1+\frac{1}{4}
h_{1}(z_1,w)\epsilon+O\big(\epsilon^2\big)\right]
\\
\phantom{T(\epsilon)}{}
=\frac{2\pi\sqrt{g\big(z^*\big)}}{\sqrt{2mgR\alpha}}\sqrt{\frac{1-z_{1}^2}{1+w^2+2wz_1}}\sqrt{\epsilon}
\left[1+O(\epsilon)\right],
\end{gather*}
where $z^*$ is some value between $z_1$ and $z_2$.
In the leading factor at $\sqrt{\epsilon}$ we have $\left|\frac{1-z_{1}^2}{1+w^2+2wz_1}\right|\leq 1$ on
any rectangle $R_{\delta}$ and the function $g(z^*)$ is decreasing for $z^*\in(-1,1)$ with supremum
$g(-1)=\frac{I_3(\alpha+1-\gamma)^2}{2(\gamma+\alpha^2-1)}$ so that
\begin{gather*}
\frac{2\pi\sqrt{g\big(z^*\big)}}{\sqrt{2mgR\alpha}}\sqrt{\frac{1-z_{1}^2}{1+w^2+2wz_1}}\sqrt{\epsilon}
<\frac{2\pi\sqrt{g(-1)}\sqrt{\epsilon}}{\sqrt{2mgR\alpha}}=2\pi\left(\frac{RI_3\gamma(\alpha+1-\gamma)}{b}
\right)=:T_{\max}.
\end{gather*}
We summarize these results in a~proposition.
\begin{proposition}\label{T_epsilon}\quad
\begin{itemize}\itemsep=0pt
\item[$i)$]For
$D\in(D_1,D_0)=\Big({-}\frac{\lambda\sqrt{d(-1)}}{R(1+\alpha)},\frac{\lambda\sqrt{d(1)}}{R(1-\alpha)}\Big)$,
$\tilde{E}>V(z_{\min},D,\lambda)$ with
\begin{gather*}
{\lambda=C\lambda_{\text{\rm thres}}
=C\frac{\sqrt{mgR^3I_3\alpha}(1+\alpha)^2}{\sqrt{1+\alpha-\gamma}}},\qquad
 C>1,
\end{gather*} the period of oscillations behaves as
\begin{gather*}
\frac{2\pi\sqrt{g(z^*)}}{\sqrt{2mgR\alpha}}\sqrt{\frac{1-z_{1}^2}{1+w^2+2wz_1}}\sqrt{\epsilon}
\left[1+O(\epsilon)\right],
\qquad
\text{as}
\quad
\epsilon=\frac{2\beta}{b^2}\to0.
\end{gather*}
\item[$ii)$] The leading factor is bounded
\[
{\frac{2\pi\sqrt{g(z^*)}}{\sqrt{2mgR\alpha}}\sqrt{\frac{1-z_{1}^2}{1+w^2+2wz_1}}
\sqrt{\epsilon}\leq\frac{2\pi RI_3\gamma(\alpha+1-\gamma)}{b}=:T_{\max}}
\]
for $w$ and $z_1$ in any closed rectangle $(z_1,w)\in[-1,1]\times[-1+\delta,1-\delta]$.
\end{itemize}
\end{proposition}

\begin{remark} Notice that the upper bound is taken uniformly w.r.t.\
both variables $(z_1,w)$ over~$R_{\delta}$.
It does not take into account the dependence of the root $z_1(a,b)$ on $a$, $b$ and that $z^*\in[z_1,z_2]$
belongs to a~subinterval of $[-1,1]$.
\end{remark}

A better estimate could be dif\/f\/icult to f\/ind due to the complexity of the expression for
$V(z,D,\lambda)$.
It is actually not needed when we do qualitative analysis of oscillations within a~deforming potential.
Here we wanted to see how the period of oscillations within the potential $V(z,D,\lambda)$ depends on the
value of Jellett's integral $\lambda$, as stated in Proposition~\ref{T_epsilon}, in order to relate the
time of inversion $T_{\text{inv}}$ to this period.
The dependence $T_{\max}\sim\frac{1}{b}\sim\frac{1}{\lambda}$ implies that the frequency of oscillations
within the potential behaves as $\frac{2\pi}{T_{\max}}\sim\lambda$.

For formulating a~suf\/f\/icient condition for having oscillating behaviour of $\theta(t)$ we need to know
that there is an upper bound for the period of oscillations within $V(z,D,\lambda)$.
To f\/ind a~universal bound independent of choice of ${\lambda>\lambda_{\text{thres}}}$ could be
dif\/f\/icult because the functions $h_1$ and $h_2$ have singularities at the boundary of the rectangle
$(z_1,w)\in[-1,1]\times[-1,1]$.
Finding a~universal bound would require detailed analysis of the interdependence between $z_1$ and $w$
during inversion.
Therefore we restrict our estimate to the region $w\in[-1+\delta,1-\delta]$ with a~certain suitable
$\delta$ and $\epsilon<\frac{1}{C^2}$, meaning $\lambda>C\lambda_{\text{thres}}$
(see~\eqref{epsilon_estimate}).

We consider period $T$ given by~\eqref{period_integral} as a~function of $\epsilon$, $w=\frac{a}{b}$ and
${z_1\in[-1,1]}$, but we drop here the assumption $\epsilon\to 0$.
Let us take $\epsilon<0.9$, ${w\in[-1+\delta,1-\delta]}$ with $\delta=0.0001$.
These are physically well justif\/ied values since $\epsilon<0.9$ means
${\lambda>1.054\lambda_{\text{thres}}}$ and $w=\pm 0.9999$ corresponds to extremely vertical initial
angular momentum $\mathbf{L}$ that is practically never taken by a~toy TT.

For $(z_1,w)\in R_{\delta}=[-1,1]\times[-1+\delta,1-\delta]$ with $\delta=0.0001$ and $\epsilon<0.9$, the
function appearing in the calculation of $\frac{1}{\sqrt{z_2-z_3}}$~\eqref{sqrtz2z3}:
\begin{gather*}
h_{3}(z_1,w)=\left(1+\frac{(1-z_{1}^2)^2}{(1+w^2+2wz_1)^2}\epsilon^2-\frac{2(1-z_{1}  %\big
^2)((1+w^2)z_1+2w)}{(1+w^2+2wz_1)^2}\epsilon\right)^{-1/4},
\end{gather*}
can be shown, by checking the maximum and minimum using a~computer algebra system since~$h_3$ does not have
singularities in $R_\delta$, to satisfy $0<h_{3}(z_1,w)<3.15$.
Then we get
\begin{gather*}
\frac{1}{\sqrt{z_2-z_3}}=\sqrt{\epsilon}\frac{1}{\sqrt{2}}\sqrt{\frac{1-z_{1}^2}{1+w^2+2wz_1}}
h_3(z_1,w)\leq2.23\sqrt{\epsilon}.
\end{gather*}
The parameter $k^2$~\eqref{k2_3} can similarly be shown to satisfy $0<k^2<0.342$ in this $R_{\delta}$,
which implies that the elliptic integral in~\eqref{complete_elliptic} satisf\/ies
$\frac{\pi}{2}<K(k^2)<1.74$.
Thus for the period $T(\epsilon)$ we have the estimate
\begin{gather*}
T(\epsilon)\leq2.23\cdot1.74\frac{4\sqrt{g\big(z^*\big)}}{\sqrt{mgR\alpha}}\sqrt{\epsilon}
<21.95\left(\frac{RI_3\gamma(\alpha+1-\gamma)}{b}\right)=:T_{\text{upp}}, %\text~152
\end{gather*}
provided that $\epsilon<0.9$.

Notice that $T_{\text{upp}}=\frac{21.95}{2\pi}T_{\max}$ is (as expected) larger than $T_{\max}$ that is
providing a~bound for the leading factor of $T(\epsilon)$.
\begin{proposition}
\label{T_upp}
Let $(z_1,w)\in R_{\delta}$ with $\delta=0.0001$, $\epsilon<0.9$ so that $b^2>\frac{20\beta}{9}$, which
corresponds to $\lambda>1.054\lambda_{\text{\rm thres}}$.
Then the period of oscillations is bounded by
$T(\epsilon)<T_{\text{\rm upp}}=21.95\left(\frac{RI_3\gamma(\alpha+1-\gamma)}{b}\right)$.
\end{proposition}

The estimate $T_{\text{upp}}$ is not the best possible estimate because it is taken uniformly over
$R_{\delta}$, but it is of the same order of magnitude as $T_{\max}$ (${T_{\text{upp}}\approx 3.5
T_{\max}}$) that also has been taken uniformly.
For large angular velocities (implying small ${\epsilon}$) $T_{\max}$ provides better estimate of magnitude
of period of oscillation.
Real oscillations within $V(z,D,\lambda)$, as observed in numerical simulations, for most of the time have
shorter period.

If the time of inversion is an order of magnitude larger than $T_{\text{upp}}$, say
${T_{\text{inv}}>10T_{\text{upp}}}$, then $\theta(t)$ changes sign many times and it has oscillatory
behaviour.

It should be stressed that behaviour of a~correctly built TT with ${1-\alpha<\gamma<1+\alpha}$ depends
strongly on the initial conditions and it also depends on the friction function
$\mu(\mathbf{s},\mathbf{L},\mathbf{\hat{3}},\mathbf{v}_A,t)$.
The majority of solutions do not demonstrate inverting behaviour and even inverting solutions do not
necessarily have to be oscillatory.
The estimate found here provide a~natural condition for the oscillatory behaviour of~TT.

Notice the dependence of $T_{\max}\sim\frac{1}{\lambda}$ and $T_{\text{upp}}\sim\frac{1}{\lambda}$.
A high initial angular momentum~$\mathbf{L}$ with $\mathbf{L}$ almost parallel to the vertical axis
$\hat{z}$ implies a~high value of Jellett's integral
$\lambda=RI_1\dot{\varphi}\sin^2\theta-RI_3\omega_3(\alpha-\cos\theta)$ since
$-RI_3\omega_3(\alpha-\cos\theta)>0$.
Then the period of nutations becomes small and $\theta(t)$ performs many oscillations within a~given
inversion time $T_{\text{inv}}$.
\begin{example}
By using again parameters from Example~\ref{main_example}, we see that $T_{\max}=\frac{2\pi R
I_3\gamma(\alpha+1-\gamma)}{b}=2.1816\cdot 10^{-8}\cdot\frac{2\pi}{b}$.
Since $b$ has the lower bound $\frac{\lambda\gamma\alpha}{1+\alpha}=1.4851\cdot 10^{-6}$ and the upper
bound $\frac{\lambda\gamma\alpha}{1-\alpha}=2.7581\cdot 10^{-6}$, the maximal period of nutation
satisf\/ies $0.0497<T_{\max}<0.0923$ seconds.
\end{example}

\section{Conclusions}

An analysis of the dynamics of an inverting top has been performed in the special case
$1-\alpha<1-\alpha^2<\gamma<1<1+\alpha$ and $I_1=\frac{I_{3}^2+mR^2I_3(1-\alpha^2)}{I_3+mR^2}$ for which
the ef\/fective potential $V(z,D(t),\lambda)$ is a~rational function of $z=\cos\theta$.

We have shown that this potential is strictly convex, has only one minimum
(Proposition~\ref{rational_convex}) and that this minimum $z_{\min}$ moves from the interval
$[1-\epsilon,1]$ to the interval $[-1,-1+\epsilon]$ (Proposition~\ref{local_min}).

The high frequency period of oscillations $T_{\max}$ and the maximal period of oscillations
$T_{\text{upp}}$ within the potential $V(z,D,\lambda)$ has been estimated in
Propositions~\ref{T_epsilon},~\ref{T_upp}.

When the time of inversion $T_{\text{inv}}$ is an order of magnitude larger than the maximal period of
oscillations $T_{\max}$ within the potential $V(z,D,\lambda)$, then the velocity $\dot{\theta}(t)$ has to
change sign many times and the solution is oscillatory.
It is visible as nutational motion of the symmetry axis $\mathbf{\hat{3}}$ on the unit sphere $S^2$ between
two latitudes while the nutational band moves from a~neighborhood of the north pole into the neighborhood
of the south pole when TT inverts.

\subsection*{Acknowledgements}

We would like to thank Hans Lundmark for discussions and for valuable help
with numerical simulations.

\pdfbookmark[1]{References}{ref}
\LastPageEnding

\end{document}